\documentclass[12pt]{amsart}
\usepackage{latexsym}
\usepackage{amsmath}
\usepackage{amssymb}
\usepackage{mathrsfs}
\input epsf
\input xy
\xyoption{all}
\input diagxy

\def\co{\operatorname{co}}
\def\Fd{\operatorname{Fd}}
\def\spt{\operatorname{spt}}
\def\Sym{\operatorname{Sym}}

\begin{document}

\title
[Multiple Criteria Problems]
{Multiple Criteria Problems\\ Over Minkowski Balls}

\author{S.~S. Kutateladze}
\address[]{
Sobolev Institute of Mathematics\newline
\indent 4 Koptyug Avenue\newline
\indent Novosibirsk, 630090\newline
\indent Russia}
\email{
sskut@math.nsc.ru
}
\begin{abstract}
Under study are some vector optimization problems over the space of Minkowski balls,
i.~e., symmetric convex compact subsets in Euclidean space.
A typical problem requires to achieve the best result in the presence of conflicting
goals; e.g., given the surface area of a symmetric convex body~$\mathfrak x$,
we try to maximize the volume of~$\mathfrak x$  and minimize the width
of~$\mathfrak x$ simultaneously.
\end{abstract}
\date{March 22, 2013}


\maketitle

\section*{Introduction}
Vector optimization is another name for multiple criteria decision making.
The mathematical technique of the field is rich but leaves much to be desired
(for instance, see~\cite{KK3}-- \cite{BGW}). One of the
reasons behind this is the fact that the classical areas of mathematics dealing
with extremal  problems pay practically no attention to the case of multiple criteria.
So it seems reasonable to suggest  attractive theoretical problems that involve
many criteria.  Some geometrical problems of the sort were considered in \cite{Kut09}.
In this article  we address similar problems over symmetric convex bodies, using the
the same technique that stems from the classical Alexandrov's approach to extremal
problems of convex geometry~\cite{SW-I}.

\section{Convex Bodies, Balls, and Dual Cones}

A~{\it convex figure\/} is a~compact convex set. A~{\it convex body\/}
is a~solid convex figure.
The {\it Minkowski  duality\/} identifies
a~convex figure $S$ in
$\mathbb R^N$ and its {\it support function\/}
$S(z):=\sup\{(x,z)\mid  x\in S\}$ for $z\in \mathbb R^N$.
Considering the members of $\mathbb R^N$ as singletons, we assume that
$\mathbb R^N$ lies in the set $\mathscr V_N$
of all compact convex subsets
of $\mathbb R^N$.

The Minkowski duality makes $\mathscr V_N$ into a~cone
in the space $C(S_{N-1})$  of continuous functions on the Euclidean unit sphere
$S_{N-1}$, the boundary of the unit ball $\mathfrak z_N$.
The
{\it linear span\/}
$[\mathscr V_N]$ of~$\mathscr V_N$ is dense in $C(S_{N-1})$, bears
a~natural structure of a~vector lattice
and is usually referred to as the {\it space of convex sets}.

The study of this space stems from the pioneering breakthrough of
Alexandrov in 1937~(see \cite{SW-I}) and the further insights of
Radstr\"{o}m, H\"{o}rmander, and Pinsker (see~\cite{KutRub}).

 A measure $\mu$ {\it linearly majorizes\/} or {\it dominates\/}
a~measure $\nu$  on $S_{N-1}$ provided that to each decomposition of
$S_{N-1}$ into finitely many disjoint Borel sets $U_1,\dots,U_m$
there are measures $\mu_1,\dots,\mu_m$ with sum $\mu$
such that every difference $\mu_k - \nu|_{U_k}$
annihilates
all restrictions to $S_{N-1}$ of linear functionals over
$\mathbb R^N$. In symbols, we write $\mu\,{\gg}{}_{\mathbb R^N} \nu$.

Reshetnyak    proved  in 1954 (see~\cite{Reshetnyak}) that
$$
\int\limits_{S_N-1} p d\mu \ge  \int\limits_{S_N-1} p d\nu
$$
for each  sublinear    functional  $p$
on  $\mathbb R^N$   if   $\mu\,{\gg}{}_{\mathbb R^N} \nu$.
This gave an important trick for generating positive linear functionals
over various classes of convex  surfaces and functions. The converse
of the Reshetnyak result was appeared in~\cite{Kut69} and \cite{Kut70}.

 Alexandrov proved the unique existence of
a translate of a convex body given its surface area function, thus
completing the solution of
the Minkowski problem.
Each surface area function is an {\it Alexandrov measure}.
So we call a positive measure on the unit sphere which is supported by
no great hypersphere and which annihilates
singletons.

Each Alexandrov measure is a translation-inva\-riant
additive functional over the cone
$\mathscr V_N$.
The cone of positive translation-invariant measures in the
dual $C'(S_{N-1})$ of
 $C(S_{N-1})$ is denoted by~$\mathscr A_N$.

 Given $\mathfrak x, \mathfrak y\in \mathscr V_N$,  the record
$\mathfrak x\,{=}{}_{\mathbb R^N}\mathfrak y$ means that $\mathfrak x$
and $\mathfrak y$  are  equal up to translation or, in other words,
are translates of one another.
So, ${=}{}_{\mathbb R^N}$ is the associate equivalence of
the preorder ${\ge}{}_{\mathbb R^N}$ on $\mathscr V_N$ of
the possibility of inserting one figure into the other
by translation.

The sum of the surface area measures of
$\mathfrak x$ and $\mathfrak y$ generates the unique class
$\mathfrak x\# \mathfrak y$ of translates which is referred to as the
{\it Blaschke sum\/} of $\mathfrak x$ and~$\mathfrak y$.
There is no need in  discriminating between a  convex figure,
the coset of its translates in  $\mathscr V_N/\mathbb R^N$,
and the corresponding measure in $\mathscr A_N$.

Let $C(S_{N-1})/\mathbb R^N$ stand for the factor space of
$C(S_{N-1})$ by the subspace of all restrictions of linear
functionals on $\mathbb R^N$ to $S_{N-1}$.
Let $[\mathscr A_N]$ be the space $\mathscr A_N-\mathscr A_N$
of translation-invariant measures, in fact, the linear span
of the set of Alexandrov measures.

$C(S_{N-1})/\mathbb R^N$ and $[\mathscr A_N]$ are made dual
by the canonical bilinear form
$$
\gathered
\langle f,\mu\rangle=\frac{1}{N}\int\limits_{S_{N-1}}fd\mu\\
(f\in C(S_{N-1})/\mathbb R^N,\ \mu \in[\mathscr A_N]).
\endgathered
$$

For $\mathfrak x\in\mathscr V_N/\mathbb R^N$ and $\mathfrak y\in\mathscr A_N$,
the quantity
$\langle {\mathfrak x},{\mathfrak y}\rangle$ coincides with the
{\it mixed volume\/}
$V_1 (\mathfrak y,\mathfrak x)$.

Consider the set ${\Sym\mathscr V}_N$ of centrally symmetric cosets of convex compact
sets. Clearly, a translation-invariant linear functional $f$ is positive over $\Sym{\mathscr V}_N$ if and only if the {\it symmetrization\/} $\Sym(f)$ is positive over $\mathscr V_N$. Here $\Sym(f)$ is the dual of the descent of the even part operator on the factor-space,  since the symmetrization of  a measure is the dual of the even part operator over $C(S_{N-1})$. We will denote the even part operator, its descent and dual by the same symbol $\Sym(\cdot)$.

 Given a cone $K$ in a vector space $X$ in duality with another vector space
$Y$,  the {\it dual\/} of $K$ is
$$
K^*:=\{y\in Y\mid (\forall x\in K)\ \langle x,y\rangle\ge 0\}.
$$

To a convex subset $U$ of $X$ and $\bar x\in U$
there corresponds
$$
U_{\bar x}:=\Fd (U,\bar x):=\{h\in X\mid (\exists \alpha \ge 0)\ \bar x+\alpha h\in U \},
$$
the {\it cone of feasible directions\/}
of $U$ at $\bar x$.

Let $\bar {\mathfrak x}\in{\mathscr A}_N$.
Then the dual  $\mathscr A^*_{N,\bar{\mathfrak x}}$ of the cone of
feasible directions of $\mathscr A_N$
at~$\bar{\mathfrak x}$ may be represented as follows
$$
{\mathscr A}^*_{N,\bar{\mathfrak x}}=\{f\in{\mathscr A}^*_N\mid
\langle\bar {\mathfrak x},f\rangle=0\}.
$$

The description of the dual of the feasible cones are  well known
(see \cite[Preposition~4.3]{Kut07}.

{\sl
Let $\mathfrak x$ and $\mathfrak y$ be  convex figures. Then

$(1)$ $\mu(\mathfrak x)- \mu(\mathfrak y)\in \mathscr V^*_N
\leftrightarrow \mu(\mathfrak x)\,{\gg}{}_{\mathbb R^N} \mu(\mathfrak y)$;

$(2)$ If $\mathfrak x\ge{}_{\mathbb R^N}\mathfrak y$
then  $\mu(\mathfrak x)\,{\gg}{}_{\mathbb R^N} \mu(\mathfrak y)$;

$(3)$ $\mathfrak x\ge{}_{\mathbb R^2}\mathfrak y\leftrightarrow
\mu(\mathfrak x)\,{\gg}{}_{\mathbb R^2} \mu(\mathfrak y)$;

$(4)$ If $\mu (\mathfrak y)-\mu (\bar{\mathfrak x})\in\mathscr V^*_{N,\bar{\mathfrak x}}$
then
$\mathfrak y=_{\mathbb R^N}\bar{\mathfrak x}$ for $\bar{\mathfrak x}\in \mathscr V_N$.
}

\noindent
From this the dual cones are available in the case of Minkowski balls.

{\sl
Let $\mathfrak x$ and $\mathfrak y$ be   convex figures. Then

$(1)$ $\mu(\mathfrak x)- \mu(\mathfrak y)\in \Sym\mathscr V^*_N
\leftrightarrow \Sym(\mu(\mathfrak x))\,{\gg}{}_{\mathbb R^N} \Sym(\mu(\mathfrak y))$;

$(2)$ If $\mathfrak x\ge{}_{\mathbb R^N}\mathfrak y$
then  $\Sym(\mu(\mathfrak x))\,{\gg}{}_{\mathbb R^N} \Sym(\mu(\mathfrak y))$;

$(3)$ $\Sym(\mathfrak x)\ge{}_{\mathbb R^2}\Sym(\mathfrak y)\leftrightarrow
\Sym(\mu(\mathfrak x))\,{\gg}{}_{\mathbb R^2} \Sym(\mu(\mathfrak y))$;

$(4)$ If $\mu (\mathfrak y)-\mu (\bar{\mathfrak x})\in(\Sym\mathscr V)^*_{N,\bar{\mathfrak x}}$
then
$\Sym(\mathfrak y)=_{\mathbb R^N}\bar{\mathfrak x}$ for $\bar{\mathfrak x}\in (\Sym\mathscr V)_N$.
}

\section{Alexandrov's Approach to the Urysohn Problem}

Alexandrov observed that the gradient of $V(\cdot)$ at $\mathfrak x$ is proportional
to $\mu(\mathfrak x)$ and so minimizing $\langle \cdot,\mu\rangle$ over $\{V=1\}$
will yield the equality $\mu=\mu(\mathfrak x)$ by the  Lagrange multiplier
rule. But this idea fails since the interior of ${\mathscr V}_{N}$ is empty.
The fact that DC-functions are dense in $C(S_{N-1})$ is not helpful at all.

Alexandrov extended the volume to the positive cone of $C(S_{N-1})$ by the formula
$V(f):=\langle f,\mu(\co(f))\rangle$ with $\co(f)$ the envelope of support functions
below $f$.  He also observed that $V(f)=V(\co(f))$. The ingenious  trick settled all for the Minkowski problem.
This was done in 1938 but still is one of the summits of convexity.

In fact, Alexandrov suggested a functional analytical approach to extremal problems
for convex surfaces. To follow it directly in the general setting  is impossible
without the above description of the dual cones. The obvious limitations of the Lagrange multiplier rule are immaterial in the case of convex programs. It should be emphasized that the classical isoperimetric problem is not a Minkowski convex program in dimensions greater than~2.
The convex counterpart is the Urysohn problem of maximizing volume given integral breadth \cite{Urysohn}.
The constraints of inclusion type are convex in the Minkowski structure,  which
opens way to complete solution  of new classes of Urysohn-type problems.

{\bf The External Urysohn Problem:}
 Among the convex figures, circumscribing $\mathfrak x_0 $ and having
integral breadth fixed, find a~convex body of greatest volume.

{\sl
A feasible convex body $\bar {\mathfrak x}$ is a solution
to~the external Urysohn problem
if and only if there are a positive  measure~$\mu $
and a positive real $\bar \alpha \in \mathbb R_+$ satisfying

$(1)$ $\bar \alpha \mu
(\mathfrak z_N)\,{\gg}{}_{\mathbb R^N}\mu (\bar {\mathfrak x})+\mu $;

$(2)$~$V(\bar {\mathfrak x})+\frac{1}{N}\int\limits_{S_{N-1}}
\bar {\mathfrak x}d\mu =\bar \alpha V_1 (\mathfrak z_N,\bar {\mathfrak x})$;

$(3)$~$\bar{\mathfrak x}(z)={\mathfrak x}_0 (z)$
for all $z$ in the support of~$\mu $, i.~e. $z\in\spt(\mu)$.
}

  If ${\mathfrak x}_0 ={\mathfrak z}_{N-1}$ then $\bar{\mathfrak x}$
is a {\it spherical lens} and $\mu$ is the restriction
of the surface area function
of the ball of radius
$\bar \alpha ^{1/(N-1)}$
to the complement of the support of the lens to~$S_{N-1}$.

If ${\mathfrak x}_0$ is an equilateral triangle then the solution
$\bar {\mathfrak x}$ looks as in Fig.~1:

$\bar {\mathfrak x}$ is the union of~${\mathfrak x}_0$
and  three congruent slices of a circle of radius~$\bar \alpha$ and
centers $O_1$--$O_3$, while
$\mu$ is the restriction of $\mu(\mathfrak z_2)$
to the subset of $S_1$ comprising the endpoints
of the unit vectors of the shaded zone.

Fig.~2 presents  the general solution of the internal Urysohn problem inside a triangle
in the class of Minkowski balls.

\centerline{\epsfxsize8cm\epsfbox{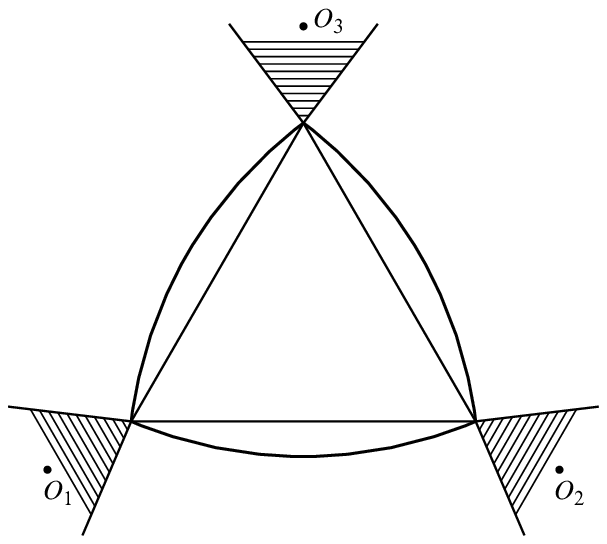}}
\centerline{Fig.~1}

\centerline{\epsfxsize6cm\epsfbox{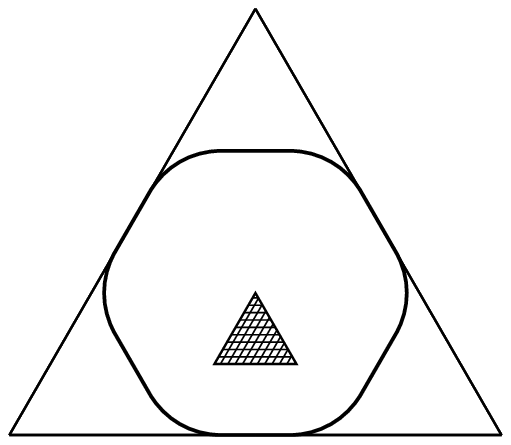}}
\centerline{Fig.~2}

\section{Pareto's Approach to Vector Optimization\\ Over Minkowski Balls}
  Consider a~bunch of  economic agents
each of which intends to maximize his own income.
The {\it Pareto efficiency principle\/}  asserts
that  as an effective agreement of the conflicting goals it is reasonable
to take any state in which nobody can increase his income in any way other
than diminishing the income of at least one of the other fellow members.
 Formally speaking, this implies the search of the maximal elements
of the set comprising the tuples of incomes of the agents
at every state; i.e., some vectors of a finite-dimensional
arithmetic space endowed with the coordinatewise order. Clearly,
the concept of Pareto optimality was already abstracted to arbitrary
ordered vector spaces.

By way of example, consider a few multiple criteria problems of isoperimetric type.
For more detail, see \cite{Kut09}.

{\bf Vector Isoperimetric Problem Over Minkowski Balls:}
 Given are some convex bodies
$\mathfrak y_1,\dots,\mathfrak y_M$.
Find a symmetric convex body $\mathfrak x$ encompassing a given volume
and minimizing each of the mixed volumes $V_1(\mathfrak x,\mathfrak y_1),\dots,V_1(\mathfrak x,\mathfrak y_M)$.
In symbols,
$$
\mathfrak x\in \Sym(\mathscr A_N);\
\widehat p(\mathfrak x)\ge \widehat p(\bar{\mathfrak x});\
(\langle\mathfrak y_1,\mathfrak x\rangle,\dots,\langle\mathfrak y_M,\mathfrak x\rangle)\rightarrow\inf\!.
$$
Clearly, this is a~Slater regular convex program in the Blaschke structure.

 {\sl
Each Pareto-optimal solution $\bar{\mathfrak x}$ of the vector isoperimetric problem
has the form}
$$
\bar{\mathfrak x}=\alpha_1\Sym({\mathfrak y}_1)+\dots+\alpha_m\Sym({\mathfrak y}_m),
$$
where $\alpha_1,\dots,\alpha_m$ are positive reals.

{\bf Internal Urysohn Problem with Flattening Over Minkowski Balls:}
Given are some~convex body
$\mathfrak x_0\in\Sym\mathscr V_N$ and some flattening direction~$\bar z\in S_{N-1}$.
Considering $\mathfrak x\subset\mathfrak x_0$ of
fixed integral breadth, maximize the volume of~$\mathfrak x$ and  minimize the
breadth of $\mathfrak x$ in the flattening direction:
$\mathfrak x\in\Sym\mathscr V_N;\
\mathfrak x\subset{\mathfrak x}_0;\
\langle \mathfrak x,{\mathfrak z}_N\rangle \ge \langle\bar{\mathfrak x},{\mathfrak z}_N\rangle;\
(-p(\mathfrak x), b_{\bar z}(\mathfrak x)) \to\inf\!.
$

{\sl For a feasible symmetric convex body $\bar{\mathfrak x}$ to be Pareto-optimal in
the internal Urysohn problem with the flattening
direction~$\bar z$  over Minkowski balls it is necessary and sufficient that there be
positive reals $\alpha$ and $\beta$ together with a~convex figure $\mathfrak x$ satisfying}
$$
\gathered
\mu(\bar{\mathfrak x})=\Sym(\mu(\mathfrak x))+ \alpha\mu({\mathfrak z}_N)+\beta(\varepsilon_{\bar z}+\varepsilon_{-\bar z});\\
\bar{\mathfrak x}(z)={\mathfrak x}_0(z)\quad (z\in\spt(\mu(\mathfrak x)).
\endgathered
$$

By way of illustration we will derive the optimality criterion
in somewhat superfluous detail.

Note firstly that the internal Urysohn problem with flattening
over Minkowski balls may be rephrased in~$C(S_{N-1})$ as the following
two-objective program
$$
\gathered
\mathfrak x\in \Sym\mathscr V_N;\\
\max\{\mathfrak x(z)- {\mathfrak x}_0(z)\mid z\in S_{N-1}\}\le 0;\\
\langle \mathfrak x,{\mathfrak z}_N\rangle \ge \langle\bar{\mathfrak x},{\mathfrak z}_N\rangle;\\
(-p(\mathfrak x), b_{\bar z}(\mathfrak x))\to\inf\!.
\endgathered
$$

The problem of Pareto optimization reduces to the
scalar program
$$
\gathered
\mathfrak x\in \Sym\mathscr V_N;\\
\max\{\max\{\mathfrak x(z)- {\mathfrak x}_0(z)\mid z\in S_{N-1}\},
\langle\bar{\mathfrak x},{\mathfrak z}_N\rangle-\langle \mathfrak x,{\mathfrak z}_N\rangle\}\le 0;\\
\max\{-p(\mathfrak x), b_{\bar z}(\mathfrak x)\}\to\inf\!.
\endgathered
$$

The last program is Slater-regular and so we may apply the
{\it Lagrange principle}. In other words, the value of the program under consideration
coincides
with the value of the free minimization problem
for an appropriate Lagrangian:
$$
\gathered
\mathfrak x\in \Sym\mathscr V_N;\\
\max\{-p(\mathfrak x), b_{\bar z}(\mathfrak x)\} +\gamma\max\{
\max\{\mathfrak x(z)- {\mathfrak x}_0(z)\mid z\in S_{N-1}\},\\
\langle\bar{\mathfrak x},{\mathfrak z}_N\rangle-\langle \mathfrak x,{\mathfrak z}_N\rangle \}\to\inf\!.
\endgathered
$$
Here $\gamma$ is a~positive Lagrange multiplier.

We are left with differentiating the Lagrangian
along the feasible directions
and appealing to~the description of the dual cones. Note in particular that
the relation
$$
\bar{\mathfrak x}(z)={\mathfrak x}_0(z)\quad (z\in\spt(\mu(\mathfrak x)))
$$
is the  {\it complementary slackness condition\/}
 standard in mathematical programming.
The proof of the optimality criterion for the Urysohn problem with
flattening over Minkowski balls complete.

{\bf Rotational Symmetry:}
Assume that a plane convex figure ${\mathfrak x}_0\in\mathscr V_2$ has the symmetry axis $A_{\bar z}$
with generator~$\bar z$.  Assume further that ${\mathfrak x}_{00}$ is the result of rotating
$\mathfrak x_0$  around the symmetry axis $A_{\bar z}$ in~$\mathbb R^3$.
Consider the problem:
$$
\gathered
\mathfrak x\in\mathscr V_3;\\
\mathfrak x  \text{\ is\ a\ convex\ body\ of\ rotation\ around}\ A_{\bar z};\\
\mathfrak x\supset{\mathfrak x}_{00};\
\langle {\mathfrak z}_N, \mathfrak x\rangle \ge \langle{\mathfrak z}_N,\bar{\mathfrak x}\rangle;\\
(-p(\mathfrak x), b_{\bar z}(\mathfrak x)) \to\inf\!.
\endgathered
$$

{\sl Each Pareto-optimal solution  is the result
of rotating around the symmetry axis a~Pareto-optimal solution of the plane internal
Urysohn problem with flattening in the direction of the axis}.

{\bf The External Urysohn Problem with Flattening Over Min\-kow\-ski Balls:}
  Given are some convex body
$\mathfrak x_0\in\mathscr V_N$ and~flattening direction~$\bar z\in S_{N-1}$.
Considering Minkowski balls $\mathfrak x\supset {\mathfrak x}_0$ of fixed integral breadth,
maximize volume and minimize  breadth in
the flattening direction:
$
\mathfrak x\in\Sym\mathscr V_N;\
\mathfrak x\supset{\mathfrak x}_0;\
\langle \mathfrak x,{\mathfrak z}_N\rangle \ge \langle\bar{\mathfrak x},{\mathfrak z}_N\rangle;\
(-p(\mathfrak x), b_{\bar z}(\mathfrak x)) \to\inf\!.
$

{\sl For a feasible convex body $\bar{\mathfrak x}$ to be a Pareto-optimal solution of
the external Urysohn problem with flattening over Minkowski balls it is necessary and
sufficient that there be
positive reals $\alpha$ and $\beta$ together with a convex figure  $\mathfrak x$ satisfying}
$$
\gathered
\mu(\bar{\mathfrak x})+\Sym(\mu(\mathfrak x))\gg{}_{{\mathbb R}^N} \alpha\mu({\mathfrak z}_N)+\beta(\varepsilon_{\bar z}+\varepsilon_{-\bar z});\\
V(\bar{\mathfrak x})+V_1(\Sym(\mathfrak x),\bar{\mathfrak x})=\alpha V_1({\mathfrak z}_N,\bar{\mathfrak x})+ 2N\beta b_{\bar z}(\bar{\mathfrak x});\\
\bar{\mathfrak x}(z)={\mathfrak x}_0(z)\quad (z\in\spt(\mu(\mathfrak x)).
\endgathered
$$

\end{document}